\documentclass[12pt,leqno,twoside]{article} 
\usepackage{bezier,amsmath,amssymb,stmaryrd,lscape}
\usepackage[colorlinks, linkcolor=black, citecolor=black, urlcolor=black, hyperfootnotes=true]{hyperref}                % hyperlinked pdf-output

\input xy
\xyoption{all}
\input{rk.sty}
\newcommand{\ulAl}{\ul{\Al}}
\newcommand{\ulmodr}{\ul{\text{\rm mod}}\text{\rm -}}
\newcommand{\veps}{\epsilon}
\newcommand{\ulEl}{\ul{\El\!}\,}

\begin{document}

\title{On adjoint functors of the Heller operator}
\author{Matthias K\"unzer}
% \date{August 31, 2007} 
\maketitle

\begin{small}
\begin{quote}
\begin{center}{\bf Abstract}\end{center}\vspace*{2mm}
Given an abelian category $\Al$ with enough projectives, we can form its stable category $\ulAl := \Al/\!\Proj(\Al)$. The Heller operator $\Og : \ulAl\lra\ulAl$ is characterised on an object $X$
by a choice of a short exact sequence $\Og X \lramono P \lraepi X$ in $\Al$ with $P$ projective. If $\Al$ is Frobenius, then $\Og$ is an equivalence, hence has a left and a right adjoint. If $\Al$
is hereditary, then $\Og$ is zero, hence has a left and a right adjoint. In general, $\Og$ is neither an equivalence nor zero. In the examples we have calculated via {\sc Magma,} it has a left adjoint, 
but in general not a right adjoint. If $\Al$ has projective covers, then $\Og$ preserves monomorphisms; this would also follow from $\Og$ having a left adjoint. I do not know an 
example where $\Og$ does not have a left adjoint.
\end{quote}
\end{small}

% this article reports on polaris.math.rwth-aachen.de:magma/sebastian/check_adjoint_1 etc.

\renewcommand{\thefootnote}{\fnsymbol{footnote}}
\renewcommand{\thefootnote}{\arabic{footnote}}

\begin{footnotesize}
\renewcommand{\baselinestretch}{0.7}%
\parskip0.0ex%
\tableofcontents%
\parskip1.2ex%
\renewcommand{\baselinestretch}{1.0}%
\end{footnotesize}%

\addtocounter{section}{-1}

\section{Introduction}
\label{SecIntro}

\subsection{The question}
\label{SecIntroQu}

Let $\El$ be an exact category. Let $\Proj(\El)$ be its full additive subcategory of relative projectives, i.e.\ for $P\in\Ob\El$ we have $P\in\Ob\Proj(\El)$ if and only
if $\liu{\El\,}{(P,-)}$ maps pure short exact sequences of $\El$ to short exact sequences of abelian groups. Suppose that $\El$ has enough relative projectives, i.e.\ suppose that for any $X\in\Ob\El$, there
exists a pure epimorphism $P\lraepi X$ in $\El$ with $P\in\Ob\Proj(\El)$.

Write $\ulEl := \El/\Proj(\El)$. The Heller operator $\Og : \ulEl\lra\ulEl$ is characterised on a given $X\in\Ob\El$ by a choice of a pure short exact sequence 
\[
\Og X\;\lramono\; P\;\lraepi\; X
\]
in $\El$ with $P$ relatively projective. This then is extended to morphisms. 

We ask whether $\Og$ has a left adjoint; cf.\ Question~\ref{QuN1}. I do not know a counterexample.

If $\El$ is a Frobenius category, then $\Og$ is an equivalence, thus has both a left and a right adjoint.

If $\El$ is hereditary, i.e.\ if $\Og \iso 0$, then $\Og$ has both a left and a right adjoint, viz.\ $0$.

\kommentar{
If $\El$ is the category of finitely generated right modules over a finite dimensional algebra $\Lm$ over a field~$k$, {\sc Auslander} and {\sc Reiten} lift $\Og$ to an exact functor $\El\lra\El$, taking the kernel
of $X\ts_k\Lm\lra X$, $x\ts\lm\lramaps x\lm\,$, which is functorial in $X\in\Ob\El\,$; cf.\ \bfcite{AR91}{Prop.\ 1.2 and \S 2}. 
}

\subsection{Monomorphisms}
\label{SecIntroMono}

If a functor has a left adjoint, then it preserves monomorphisms. So first of all, we ask whether the Heller operator $\Og: \ulEl\lra\ulEl$ preserves monomorphisms.

It turns out that if $\El$ is weakly idempotent complete and has relative projective covers in the sense of \S\ref{SecProjCov}, then $\Og$ maps monomorphisms even to coretractions; cf.\ Proposition~\ref{PropM2}.

\subsection{Construction of a left adjoint to the Heller operator $\Og$}
\label{SecIntroConstr}

Let $p\in [2,997]$ be a prime. Let $R := \Fu{p}[X]$ and $\pi := X$. 

Let $A := (R/\pi^3)(e\lraa{a} f) \iso \smatzz{R/\pi^3}{R/\pi^3}{0}{R/\pi^3}$. Let $\El := \modr A\,$.

Using {\sc Magma}~\bfcit{CaEtAl}, we construct a left adjoint $S : \ulEl\lra\ulEl$ to $\Og$. We do so likewise for certain factor rings of~$A\,$.
Cf.\ Propositions~\ref{PropE2},~\ref{PropE4}~and~\ref{PropE6}.

Let now $k$ be a field, $R := k[X]$ and $\pi := X$. Let $n\ge 1$. An $(R/\pi^n)(e\lraa{a} f)$-module is given by a morphism $X\lraa{f} Y$ in $\modr(R/\pi^n)$. 
The full subcategory of $\modr(R/\pi^n)(e\lraa{a} f)$ consisting of injective morphisms $X\lraa{f} Y$ as modules has been intensely studied; it is of finite type if $n\le 5$, 
tame if $n = 6$, wild if $n\ge 7$; cf.\ \bfcite{SR08}{(0.1), (0.6)}.

\subsection{Two counterexamples}
\label{SecIntroCounter}

The functor $\Og : \ulEl\lra\ulEl$ does not have a right adjoint in general; cf.\ Remark~\ref{RemE8}.

If existent, the functor $\Og\0 S : \ulEl\lra\ulEl$ is not idempotent in general; cf.\ Remark~\ref{RemE7}.

\subsection{Acknowledgements}
\label{SecAck}

{\sc Sebastian Thomas} asked whether there exists a category with set of weak equivalences $(\Cl,W)$ that carries the structure of a Brown fibration category, but whose $\Og$
on $\Cl[W^{-1}]$~(\footnote{Cf.\ also~\bfcite{Se03}{p.\ 210}.}) does not have a left adjoint. In our exact category context, this is Question~\ref{QuN1}. 

I thank {\sc Steffen K\"onig} for help with monomorphisms, cf.\ \S\ref{SecMonoPres}. 
I thank {\sc Sebastian Thomas} and {\sc Markus Kirschmer} for help with {\sc Magma.} I thank {\sc Markus Schmidmeier} for help with $\modr(R/\pi^n)(e\lraa{a} f)$.

\subsection{Notations and conventions}
\label{SecNotConv}

\begin{footnotesize}
\begin{itemize}
\item Given $a,\,b\,\in\,\Z$, we write $[a,b] := \{\,z\in\Z\, :\,a\le z\le b\,\}$.
\item Composition of morphisms is written naturally, $(\lraa{a}\lraa{b}) = \lraa{ab} = \lraa{a\cdot b}$. 

Composition of functors is written traditionally, $(\lraa{F}\lraa{G}) = \lraa{G\0 F}$.
\item In a category $\Cl$, given $X,\,Y\,\in\,\Ob\Cl$, we write $\liu{\Cl}{(X,Y)}$ for the set of morphisms from $X$ to $Y$.
\item Given an isomorphism $f$, we write $f^-$ for its inverse.
\item In an additive category, a morphism of the form $X\lrafl{23}{\smatez{1}{0}} X\ds Y$, or isomorphic to such a morphism, is called split monomorphic; a morphism of the form $X\ds Y\lrafl{32}{\smatze{1}{0}} X\ru{6.5}$, 
or isomorphic to such a morphism, is called split epimorphic.
\item In exact categories, pure monomorphisms are denoted by $\lramono$, pure epimorphisms by $\lraepi$ and pure squares, i.e.\ bicartesian squares with pure short exact diagonal sequence, by a box $\Box$ in the diagram.
\item Given a ring $A$, an {\it $A$-module} is a finitely generated right $A$-module.
\item Given a commutative ring $A$ and $a\in A$, we often write $A/a := A/(a) = A/aA$.
\item Given a noetherian ring $A$, we write $\ulmodr A := \ul{\modr A}$ for the factor category of $\modr A$ modulo the full additive subcategory of projectives. 

So in the language of \S\ref{SecMonoNot} below, 
we consider the abelian category $\modr A$ as an exact category with all short exact sequences declared to be pure and write $\ulmodr A$ for its classical stable category.
\end{itemize}
\end{footnotesize}
    % Introduction

\section{The Heller operator $\Omega$}
\label{SecMono}

\subsection{Notation}
\label{SecMonoNot}

Let $\El$ be an exact category in the sense of {\sc Quillen} \bfcite{Qu73}{p.\ 99} with enough relative projectives. We will use the notation of \bfcite{Ku05}{\S A.2} concerning pure short exact sequences, pure 
monomorphisms and pure epimorphisms.

Let $\Proj(\El)\tm\El$ denote the full subcategory of relative projectives. Let 
\[
\ulEl \; := \; \El/\Proj(\El)
\]
denote the classical stable category of $\El$. The residue class functor shall be denoted by
\[
\barcl
\El           & \lra      & \ulEl                    \\
(X\lraa{f} Y) & \lramaps  & (X\lrafl{25}{[f]} Y)\; . \\
\ea
\]

For each $X\in\Ob\El$, we choose a pure short exact sequence
\hypertarget{ast.1}{}
\[
\Og X\;\lramonoa{i_X}\; \PP X \;\lraepifl{34}{p_X}\; X
\leqno (\ast)
\]
with $\PP X$ relatively projective. Let the {\it Heller operator} \bfcit{He60}
\[
\Og \;:\; \ulEl \;\lra\; \ulEl
\]
be defined on the objects by the choice just made. Suppose given a morphism $X\lrafl{25}{[f]} Y$ in $\ulEl$. Choose a morphism 
\[
\xymatrix{
\Og X\arm[r]^*+<0mm,1.5mm>{\scm i_X}\ar[d]^{f'} & \PP X\are[r]^*+<0mm,1.5mm>{\scm p_X}\ar[d]^{\h f} & X\ar[d]^f \\
\Og Y\arm[r]^*+<0mm,1.5mm>{\scm i_Y}            & \PP Y\are[r]^*+<0mm,1.5mm>{\scm p_Y}              & Y         \\
}
\]
of pure short exact sequences. Let
\[
\Og[f] \; := \; [f']\; .
\]

Different choices of pure short exact sequences \hyperlink{ast.1}{$(\ast)$} yield mutually isomorphic Heller operators.

\begin{Question}
\label{QuN1}
Does $\,\Og$ have a left adjoint?
\end{Question}

I do not know a counterexample.

\subsection{Preservation of monomorphisms}
\label{SecMonoPres}

If $\Og : \ulEl \lra \ulEl$ has a left adjoint, then it preserves monomorphisms. So if, for some $\El$, the functor $\Og$ did not preserve monomorphisms, then $\Og$ could not have a left adjoint.
Under certain finiteness assumptions, however, we will show that $\Og$ maps monomorphisms to coretractions, so in particular to monomorphisms. This is to be compared to the case of $\El$ being Frobenius, where in the triangulated
category $\ulEl$ all monomorphisms are split.

\begin{Lemma}
\label{LemM1}
Suppose that for $X\in\Ob\El$ and for $s\in\liu{\El}{(\PP X,\PP X)}$ such that $sp_X = p_X\,$, the endomorphism $s$ is an isomorphism.

Suppose given $X\lraa{f} Y$ in $\El$.

If $\,[f]$ is a monomorphism, then $\Og[f]$ is a coretraction.

In particular, $\Omega$ preserves monomorphisms.
\end{Lemma}

{\it Proof.} Choose a morphism of pure short exact sequences as shown below. Insert a pullback $(T,X,\,\PP Y,Y)$ and the induced morphism $\PP X\lraa{v} T$, having $vg = \h f$ and $vq = p_X\,$. Insert a
kernel $j$ of $q$ with $jg = i_Y\,$.
\[
\xymatrix{
\Og X\arm[rrr]^*+<0mm,1.5mm>{\scm i_X}\ar[ddd]^{f'}  &&& \PP X\are[rrr]^*+<0mm,1.5mm>{\scm p_X}\ar[ddd]^(0.3){\h f}|(0.54)\hole\ar[dr]^v &                                          && X\ar[ddd]^f \\
                                                     &&&                                                                                 & T\are[urr]^*+<3mm,0mm>{\scm q}\ar[ddl]_g &&             \\
                                                     &&&                                                                                 &                                          &&             \\
\Og Y\arm[rrr]^*+<0mm,1.5mm>{\scm i_Y}\arm[uurrrr]^j &&& \PP Y\are[rrr]^*+<0mm,1.5mm>{\scm p_Y}                                          &                                          && Y           \\
}
\]
We have $f'j = i_X v$, since $f'j q = 0 = i_X p_X = i_X vq$ and $f'jg = f'i_Y = i_X\h f = i_X vg$.

We have $[q][f] = [gp_Y] = 0$. Since $[f]$ is monomorphic, we infer that $[q] = 0$. Hence there exists $T\lraa{u}\PP X$ such that $u p_X = q$. On the kernels, we obtain $\Og Y\lraa{u'}\Og X$ such that
$u' i_X = ju$.
\[
\xymatrix{
\Og X\arm[rrr]^*+<0mm,1.5mm>{\scm i_X}\ar[ddd]^{f'}                      &&& \PP X\are[rrr]^*+<0mm,1.5mm>{\scm p_X}\ar[ddd]^(0.3){\h f}|(0.54)\hole\ar[dr]^v &                                                                && X\ar[ddd]^f \\
                                                                         &&&                                                                                 & T\are[urr]^*+<3mm,0mm>{\scm q}\ar[ddl]_g\ar@/^3mm/[ul]^(0.35)u &&             \\
                                                                         &&&                                                                                 &                                                                &&             \\
\Og Y\arm[rrr]^*+<0mm,1.5mm>{\scm i_Y}\arm[uurrrr]^j\ar@/^3mm/[uuu]^{u'} &&& \PP Y\are[rrr]^*+<0mm,1.5mm>{\scm p_Y}                                          &                                                                && Y           \\
}
\]
We have $vup_X = vq = p_X\,$. Hence $vu$ is an isomorphism by assumption.

We obtain $f'u'i_X = f'ju = i_X vu$. So $(f'u',vu,\id_X)$ is a morphism of pure short exact sequences. Hence $f'u'$ is an isomorphism. Thus $f'$ is a coretraction. We conclude that $\Og [f] = [f']$ is a coretraction.
\qed

\subsection{Relative projective covers}
\label{SecProjCov}

Suppose $\El$ to be weakly idempotent complete; cf.\ \bfcite{Bu10}{Def.\ 7.2}.

A morphism $S\lraa{i} M$ in $\El$ is called {\it small} if in each pure square
\[
\xymatrix{
A\ar[r]\ar[d]\ar@{}[dr]|\Box & T\ar[d]^t \\
S\ar[r]^i                    & M         \\
}
\]
in $\El$, the morphism $T\lraa{t} M$ is purely epimorphic; cf.\ \bfcite{Ro91}{Def.\ 2.8.30}. In other words, $i$ is small iff $\smatze{i}{t}$ being purely epimorphic entails $t$ being purely epimorphic for each morphism $t$ with 
same target as~$i$. E.g.\ $i = 0$ is small, for $\smatze{0}{t} = \smatze{0}{1}t$ is purely epimorphic only if $t$ is; cf.~\bfcite{Bu10}{Prop.\ 2.16}. 

If $S\lramonofl{25}{i} M$ is small and split monomorphic, then there exists
\[
\xymatrix{
0\ar[r]\ar[d]\ar@{}[dr]|\Box & S'\arm[d]^*+<1mm,0mm>{\scm i'} \\
S\arm[r]^*+<0mm,1mm>{\scm i} & M\zw{,}                        \\
}
\]
forcing $i'$ to be an isomorphism and thus $S$ to be isomorphic to $0$.

Given $\w S\lraa{j} S\lraa{i} M$ with $S\lraa{i} M$ small, then $\w S\lraa{ji} M$ is small. In fact, given $t$ such that $\smatze{ji}{t}$ is a pure epimorphism,
the factorisation $\smatze{ji}{t} = \smatzz{j}{0}{0}{1}\smatze{i}{t}$ shows that $\smatze{i}{t}$ is a pure epimorphism; cf.\ \bfcite{Bu10}{Cor.\ 7.7}. Thus $t$ is purely epimorphic by smallness of $i$.

A {\it relative projective cover} of $X\in \Ob\El$ is a pure epimorphism $P\lraepia{p} X$ in $\El$ such that $P$ is relatively projective and such that $\Kern p \lramono P$ is small; cf.\ \bfcite{Ro91}{2.8.31}. 

We say that $\El$ {\it has relative projective covers} if for each $X\in\Ob\El$, there exists a relative projective cover $P\lraepia{p} X$.

\begin{Lemma}
\label{LemM1_5}
Suppose given a relative projective cover $P\lraepia{p} X$ in $\El$. Suppose given $P\lraa{s} P$ such that $sp = p$. Then $s$ is an isomorphism.
\end{Lemma}

{\it Proof.} We complete to a pure short exact sequence $K\lramonoa{k} P\lraepia{p} X\ru{6}$. We obtain a morphism
\[
\xymatrix{
K\arm[r]^*+<0mm,1mm>{\scm k}\ar[d]_g\ar@{}[dr]|\Box  & P\are[r]^*+<0mm,1mm>{\scm p}\ar[d]_s & X\ar@{=}[d] \\
K\arm[r]^*+<0mm,1mm>{\scm k}                         & P\are[r]^*+<0mm,1mm>{\scm p}         & X\;\zw{.}   \\
}
\]
of pure short exact sequences. Since the left hand side quadrangle is a pure square, we conclude that $s$ is purely epimorphic by smallness of $\smash{K\lramonofl{23}{k} P}$.  
Hence $s$ is split epimorphic by relative projectivity of $P\,$; cf.\ \bfcite{Bu10}{Rem.\ 7.4}. Let $\smash{L\lramonofl{23}{\ell} P}$ be a kernel of $s$. Since $\ell$ factors over the small morphism $k$, it is small as well. 
Since $\ell$ is split monomorphic, we have $L\iso 0$. Thus $s$ is an isomorphism.\qed

\begin{Proposition}
\label{PropM2}
Suppose that the exact category $\El$ is weakly idempotent complete and has relative projective covers. 

Then $\Og : \ulEl \lra \ulEl$ maps each monomorphism to a coretraction. In particular, $\Omega$ preserves monomorphisms.
\end{Proposition}

{\it Proof.} We may use relative projective covers to construct $\Og$ in \hyperlink{ast.1}{$(\ast)$}. Then Lemma~\ref{LemM1_5} allows us to apply Lemma~\ref{LemM1}.\qed
     % Monomorphisms

\section{Examples for adjoints of the Heller operator $\Og$}
\label{SecEx}

Let $R$ be a principal ideal domain, with a maximal ideal generated by an element $\pi\in R$.

Let 
\[
A \; :=\; (R/\pi^3)(e\lraa{a} f)\; .
\]
I.e.\ $A$ is the path algebra of $e\lraa{a} f$ over the ground ring $R/\pi^3$. It has primitive idempotents $e$ and $f$, and $a\in e A f$.

An object in $\modr A$ is given by a morphism $X\lra Y$ in $\modr(R/\pi^3)$. A morphism in $\modr A$ is given by a commutative quadrangle in $\modr(R/\pi^3)$.

\subsection{Example of a left adjoint}
\label{SecExLeftAdjA}

\subsubsection{A list of indecomposables}
\label{SecExListA}

Define the following objects in $\modr A$.
% Data: cf. polaris.math.rwth-aachen.de:magma/sebastian/check_adjoint_1 resp. /check_adjoint_1_loop for the epsilon values further below
\[
\ba{lclclcl}
P_1    & := & (R/\pi^3 \lrafl{25}{1} R/\pi^3)                                                      && P_2    & := & (0 \lrafl{25}{} R/\pi^3)                                                             \vsp{5}\\
X_1    & := & (R/\pi\lrafl{25}{1} R/\pi)                                                           && X_{14} & := & (R/\pi^3\lrafl{25}{\pi^2} R/\pi^3)                                                   \vsp{2}\\
X_2    & := & (R/\pi^2\lrafl{25}{1} R/\pi^2)                                                       && X_{15} & := & (R/\pi^3\lrafl{25}{} 0)                                                              \vsp{2}\\
X_3    & := & (R/\pi^2\lrafl{25}{1} R/\pi)                                                         && X_{16} & := & (R/\pi^2\ds R/\pi^3\lrafl{36}{\smatzz{1}{\pi}{1}{0}} R/\pi\ds R/\pi^3)               \vsp{2}\\
X_4    & := & (R/\pi^3\lrafl{25}{1} R/\pi^2)                                                       && X_{17} & := & (R/\pi\ds R/\pi^3\;\lrafl{36}{\smatzz{\pi}{\smash{\pi^2}}{1}{0}} R/\pi^2\ds R/\pi^3) \vsp{2}\\
X_5    & := & (R/\pi^2\lrafl{25}{\pi} R/\pi^3)                                                     && X_{18} & := & (R/\pi\ds R/\pi^3\;\lrafl{36}{\smatzz{0}{\smash{\pi^2}}{1}{\pi}} R/\pi\ds R/\pi^3)   \vsp{2}\\
X_6    & := & (R/\pi\lrafl{25}{\pi} R/\pi^2)                                                       && X_{19} & := & (R/\pi\lrafl{25}{\pi^2} R/\pi^3)                                                     \vsp{2}\\
X_7    & := & (R/\pi^2\lrafl{25}{\smatez{1}{\pi}} R/\pi\ds R/\pi^3)                                && X_{20} & := & (R/\pi^3\lrafl{25}{1} R/\pi)                                                         \vsp{2}\\
X_8    & := & (R/\pi\ds R/\pi^3\lrafl{36}{\smatze{\pi}{1}} R/\pi^2)                                && X_{21} & := & (R/\pi\lrafl{25}{} 0)                                                                \vsp{2}\\
X_9    & := & (R/\pi^3\lrafl{25}{\pi} R/\pi^3)                                                     && X_{22} & := & (0\lrafl{25}{} R/\pi)                                                                \vsp{2}\\
X_{10} & := & (R/\pi^2\lrafl{25}{} 0)                                                              && X_{23} & := & (R/\pi^3\lrafl{25}{\pi} R/\pi^2)                                                     \vsp{2}\\
X_{11} & := & (0\lrafl{25}{} R/\pi^2)                                                              && X_{24} & := & (R/\pi^2\lrafl{25}{\pi^2} R/\pi^3)                                                   \vsp{2}\\
X_{12} & := & (R/\pi\ds R/\pi^3\lrafl{36}{\smatze{\smash{\pi^2}}{\pi}} R/\pi^3)                    && X_{25} & := & (R/\pi^2\lrafl{25}{\pi} R/\pi^2)                                                     \vsp{2}\\
X_{13} & := & (R/\pi^3\lrafl{28}{\smatez{1}{\pi}} R/\pi\ds R/\pi^3)                                &&        &    &                                                                                      \vsp{2}\\
\ea
\]

A matrix inspection yields the  

\begin{Lemma}
\label{LemE1}\Absit
\begin{itemize}
\item[\rm (1)] For each projective indecomposable $A$-module $P$, there exists a unique $i\in [1,2]$ such that \mb{$P\iso P_i\,$}.
\item[\rm (2)] For each nonprojective indecomposable $A$-module $X$, there exists a unique $i\in [1,25]$ such that \mb{$X\iso X_i\,$}.
\end{itemize}
\end{Lemma}

\subsubsection{Construction of a left adjoint}
\label{SecConstrLeftAdjA}

Our aim in this section is to computationally verify the

\begin{Proposition}
\label{PropE2}
Suppose given a prime $p\in [2,997]$. Suppose that $R = \Fu{p}[X]$ and $\pi = X$.

Then the Heller operator $\,\Og : \ulmodr A\lra\ulmodr A$ has a left adjoint.
\end{Proposition}

For ease of {\sc Magma} input, we have used that 
\[
A \;\iso\; \Fu{p}\!\left(
\xymatrix{
e\ar[r]^a\ar@<0mm>@(ul,dl)[]_{u} & f\ar@<0mm>@(ur,dr)[]^{v} \\
}
\right)/ (u^3,\;v^3,\;ua - av)  
\]
as $\Fu{p}$-algebras. 

To reduce the calculation of this adjoint functor to the proof of the representability of certain functors, we use 

\begin{Lemma}[\bfcite{Sc72}{16.4.5, 4.5.1}]
\label{LemE3}
Suppose given categories $\Cl$ and $\Dl$ and a functor $\,\Cl\lraa{F}\Dl$. 

Suppose that
\[
\barcl
\Cl            & \mrafl{30}{\liu{\Dl}{(Y,F(-))}\;} & \Sets                                                                                  \\
(X\lraa{x} X') & \mramaps                          & \smash{\big(\liu{\Dl}{(Y,FX)} \mrafl{30}{\liu{\Dl}{(Y,\,Fx)}} \liu{\Dl}{(Y,FX')}\big)} \\
\ea
\]
is representable for each $Y\in\Ob\Dl$. 

Then $F$ has a left adjoint.

More precisely, given a map $\Ob\Cl\llaa{\ga}\Ob\Dl$ and an isomorphism 
\[
\liu{\Dl}{(Y,F(-))} \;\lraisoa{\phi_Y}\; \liu{\Cl}{(Y\ga,-)}
\]
for $Y\in\Ob\Dl$, there exists a left adjoint $\Cl\llaa{G}\Dl$ to $F$, i.e.\ $G\adj F$, such that, writing
\[
\eps Y \; :=\; (1_{Y\ga})(\ph_Y(Y\ga))^- \; :\; Y \; \lra \; F(Y\ga) 
\]
for $Y\in\Ob\Dl$, we have 
\[
G(Y\xrightarrow{y} Y') \;\=\; (Y\ga\xrightarrow{(y\,\cdot\,\eps Y') (\ph_Y(Y'\ga))}Y'\ga) 
\]
for $Y\lraa{y} Y$ in $\Dl$. 
\end{Lemma} 

Thus in order to construct the left adjoint to $\Og$ on $\ulmodr A$, it suffices to show that the functor $\liu{\ulmodr A\,}{(X_i\,,\,\Og(-))}$ is 
representable $i\in [1,25]$. We shall do so by an actual construction of an isotransformation from a Hom-functor. 

Suppose given $i\in [1,25]$. Such an isotransformation is necessarily of the form
\[
\barcl
\liu{\ulmodr A\,}{(S X_i\,, Y)} & \lra     & \liu{\ulmodr A\,}{(X_i\,, \Og Y)} \\
{[f]}                           & \lramaps & [\veps_i]\cdot \Og {[f]}          \\
\ea
\]
for some $SX_i\in\Ob\ulmodr A$ and some $A$-linear map $\veps_i : X_i \lra \Og S X_i\,$, where $Y\in\Ob\ulmodr A$.

So it suffices to find an $A$-module $SX_i$ and an $A$-linear map $\veps_i : X_i \lra \Og S X_i$ such that the induced map
\hypertarget{ast.2}{}
\[
\barcl
\liu{\ulmodr A\,}{(S X_i\,, X_j)} & \lra     & \liu{\ulmodr A\,}{(X_i\,, \Og X_j)} \\
{[f]}                             & \lramaps & [\veps_i]\cdot \Og {[f]}            \\
\ea
\leqno (\ast\ast)
\]
is an isomorphism for $j\in [1,25]$.

In particular, given an automorphism $\al$ of $X_i$ in $\ulmodr A\,$, an automorphism $\be$ of $SX_i$ in $\ulmodr A$ and a valid such morphism $[\veps_i]\,$, then $\al\cdot [\veps_i]\cdot \Og\be$ is another valid such morphism, sometimes
of a simpler shape.

To show that a guess for $S X_i$ is in fact the sought-for representing object, we make use of the fact that $\liu{\ulmodr A\,}{(X_i\,, \Og SX_i)}$ is finite, so that we have only a finite set of 
candidates for $\veps_i\,$. Then to check whether the candidate-induced maps \hyperlink{ast.2}{$(\ast\ast)$} are isomorphisms, is also feasible via {\sc Magma,} 
using in particular its commands \verb|ProjectiveCover|, \verb|AHom| and \verb|PHom|; cf.\ \bfcit{CaEtAl}.

We obtain
\[
\ba{lclclclclclclcl}
SX_1    & = & X_2    & \hsp{3} & \Og SX_1    & = & X_1    & \hsp{5} & SX_{14} & = & X_5           & \hsp{3} & \Og SX_{14} & = & X_{19}          \\
SX_2    & = & X_1    &         & \Og SX_2    & = & X_2    &         & SX_{15} & = & 0             &         & \Og SX_{15} & = & 0               \\
SX_3    & = & X_2    &         & \Og SX_3    & = & X_1    &         & SX_{16} & = & X_2\ds X_{19} &         & \Og SX_{16} & = & X_1\ds X_5      \\
SX_4    & = & X_1    &         & \Og SX_4    & = & X_2    &         & SX_{17} & = & X_1\ds X_5    &         & \Og SX_{17} & = & X_2\ds X_{19}   \\
SX_5    & = & X_{19} &         & \Og SX_5    & = & X_5    &         & SX_{18} & = & X_6           &         & \Og SX_{18} & = & X_7             \\
SX_6    & = & X_7    &         & \Og SX_6    & = & X_6    &         & SX_{19} & = & X_5           &         & \Og SX_{19} & = & X_{19}          \\
SX_7    & = & X_6    &         & \Og SX_7    & = & X_7    &         & SX_{20} & = & X_2           &         & \Og SX_{20} & = & X_1             \\
SX_8    & = & X_1    &         & \Og SX_8    & = & X_2    &         & SX_{21} & = & 0             &         & \Og SX_{21} & = & 0               \\
SX_9    & = & X_{19} &         & \Og SX_9    & = & X_5    &         & SX_{22} & = & X_{11}        &         & \Og SX_{22} & = & X_{22}          \\
SX_{10} & = & 0      &         & \Og SX_{10} & = & 0      &         & SX_{23} & = & X_7           &         & \Og SX_{23} & = & X_6             \\
SX_{11} & = & X_{22} &         & \Og SX_{11} & = & X_{11} &         & SX_{24} & = & X_5           &         & \Og SX_{24} & = & X_{19}          \\
SX_{12} & = & X_{19} &         & \Og SX_{12} & = & X_5    &         & SX_{25} & = & X_7           &         & \Og SX_{25} & = & X_6             \\
SX_{13} & = & X_6    &         & \Og SX_{13} & = & X_7    &         &         &   &               &         &             &   &                 \\
\ea
\]
and
\begin{footnotesize}
\[
\ba{lllllll}
\left(\raisebox{6mm}{\xymatrix@R-4mm{X_1\ar[d]_{\veps_1}\\\Og SX_1}}\right) 
& = & 
\left(\raisebox{6mm}{\xymatrix@R-4mm{
R/\pi\ar[r]^1\ar[d]_1 & R/\pi\ar[d]^1 \\
R/\pi\ar[r]^1         & R/\pi         \\
}}\right)
& \hsp{5} &
\left(\raisebox{6mm}{\xymatrix@R-4mm{X_{14}\ar[d]_{\veps_{14}}\\\Og SX_{14}}}\right) 
& = & 
\left(\raisebox{6mm}{\xymatrix@R-4mm{
R/\pi^3\ar[r]^{\pi^2}\ar[d]_1 & R/\pi^3\ar[d]^1 \\
R/\pi\ar[r]^{\pi^2}           & R/\pi^3         \\
}}\right) \\
\left(\raisebox{6mm}{\xymatrix@R-4mm{X_2\ar[d]_{\veps_2}\\\Og SX_2}}\right) 
& = & 
\left(\raisebox{6mm}{\xymatrix@R-4mm{
R/\pi^2\ar[r]^1\ar[d]_1 & R/\pi^2\ar[d]^1 \\
R/\pi^2\ar[r]^1         & R/\pi^2         \\
}}\right)
& &
\left(\raisebox{6mm}{\xymatrix@R-4mm{X_{15}\ar[d]_{\veps_{15}}\\\Og SX_{15}}}\right) 
& = & 
\left(\raisebox{6mm}{\xymatrix@R-4mm{
R/\pi^3\ar[r]\ar[d] & 0\ar[d] \\
0\ar[r]             & 0       \\
}}\right) \\
\left(\raisebox{6mm}{\xymatrix@R-4mm{X_3\ar[d]_{\veps_3}\\\Og SX_3}}\right) 
& = & 
\left(\raisebox{6mm}{\xymatrix@R-4mm{
R/\pi^2\ar[r]^1\ar[d]_1 & R/\pi\ar[d]^1 \\
R/\pi\ar[r]^1           & R/\pi         \\
}}\right)
& &
\hsp{-12}\left(\raisebox{6mm}{\xymatrix@R-4mm{X_{16}\ar[d]_{\veps_{16}}\\\Og SX_{16}}}\right) 
& \hsp{-12}= & 
\hsp{-12}\left(\raisebox{6mm}{\xymatrix@R-4mm{
R/\pi^2\ds R/\pi^3\ar[r]^{\smatzz{1}{\pi}{1}{0}}\ar[d]_{\rsmatzz{1}{1}{1}{0}} & R/\pi\ds R/\pi^3\ar[d]^{\smatzz{1}{0}{0}{1}} \\
R/\pi\ds R/\pi^2\ar[r]^{\smatzz{1}{0}{0}{\pi}}                                & R/\pi\ds R/\pi^3                             \\
}}\right) \\
\ea
\]
\[
\ba{lllllll}
\left(\raisebox{6mm}{\xymatrix@R-4mm{X_4\ar[d]_{\veps_4}\\\Og SX_4}}\right) 
& = & 
\left(\raisebox{6mm}{\xymatrix@R-4mm{
R/\pi^3\ar[r]^1\ar[d]_1 & R/\pi^2\ar[d]^1 \\
R/\pi^2\ar[r]^1         & R/\pi^2         \\
}}\right)
& &
\hsp{-12}\left(\raisebox{6mm}{\xymatrix@R-4mm{X_{17}\ar[d]_{\veps_{17}}\\\Og SX_{17}}}\right) 
& \hsp{-12}= & 
\hsp{-12}\left(\raisebox{6mm}{\xymatrix@R-4mm{
R/\pi\ds R/\pi^3\ar[r]^{\smatzz{\pi}{\;\pi^2}{1}{0}}\ar[d]_{\rsmatzz{\pi}{1}{1}{0}} & R/\pi^2\ds R/\pi^3\ar[d]^{\smatzz{1}{0}{0}{1}} \\
R/\pi^2\ds R/\pi\ar[r]^{\smatzz{1}{0}{0}{\;\pi^2}}                                  & R/\pi^2\ds R/\pi^3                             \\
}}\right) \vsp{1}\\
\left(\raisebox{6mm}{\xymatrix@R-4mm{X_5\ar[d]_{\veps_5}\\\Og SX_5}}\right) 
& = & 
\left(\raisebox{6mm}{\xymatrix@R-4mm{
R/\pi^2\ar[r]^\pi\ar[d]_1 & R/\pi^3\ar[d]^1 \\
R/\pi^2\ar[r]^\pi         & R/\pi^3         \\
}}\right)
& &
\hsp{-12}\left(\raisebox{6mm}{\xymatrix@R-4mm{X_{18}\ar[d]_{\veps_{18}}\\\Og SX_{18}}}\right) 
& \hsp{-12}= & 
\hsp{-12}\left(\raisebox{6mm}{\xymatrix@R-4mm{
R/\pi\ds R/\pi^3\ar[r]^{\smatzz{0}{\;\pi^2}{1}{\pi}}\ar[d]_{\smatze{\pi}{1}} & R/\pi\ds R/\pi^3\ar[d]^{\smatzz{1}{0}{0}{1}} \\
R/\pi^2\ar[r]^-{\smatez{1}{\pi}}                                             & R/\pi\ds R/\pi^3                             \\
}}\right) \\
\left(\raisebox{6mm}{\xymatrix@R-4mm{X_6\ar[d]_{\veps_6}\\\Og SX_6}}\right) 
& = & 
\left(\raisebox{6mm}{\xymatrix@R-4mm{
R/\pi\ar[r]^\pi\ar[d]_1 & R/\pi^2\ar[d]^1 \\
R/\pi\ar[r]^\pi         & R/\pi^2         \\
}}\right)
& &
\left(\raisebox{6mm}{\xymatrix@R-4mm{X_{19}\ar[d]_{\veps_{19}}\\\Og SX_{19}}}\right) 
& = & 
\left(\raisebox{6mm}{\xymatrix@R-4mm{
R/\pi\ar[r]^{\pi^2}\ar[d]_1 & R/\pi^3\ar[d]^1 \\
R/\pi\ar[r]^{\pi^2}         & R/\pi^3         \\
}}\right)\vsp{1} \\
\left(\raisebox{6mm}{\xymatrix@R-4mm{X_7\ar[d]_{\veps_7}\\\Og SX_7}}\right) 
& = & 
\left(\raisebox{6mm}{\xymatrix@R-4mm{
R/\pi^2\ar[r]^-{\smatez{1}{\pi}}\ar[d]_1 & R/\pi\ds R/\pi^3\ar[d]^{\smatzz{1}{0}{0}{1}} \\
R/\pi^2\ar[r]^-{\smatez{1}{\pi}}         & R/\pi\ds R/\pi^3                             \\
}}\right)
& & 
\left(\raisebox{6mm}{\xymatrix@R-4mm{X_{20}\ar[d]_{\veps_{20}}\\\Og SX_{20}}}\right) 
& = & 
\left(\raisebox{6mm}{\xymatrix@R-4mm{
R/\pi^3\ar[r]^1\ar[d]_1 & R/\pi\ar[d]^1 \\
R/\pi\ar[r]^1           & R/\pi         \\
}}\right) \vsp{1}\\
\left(\raisebox{6mm}{\xymatrix@R-4mm{X_8\ar[d]_{\veps_8}\\\Og SX_8}}\right) 
& = & 
\left(\raisebox{6mm}{\xymatrix@R-4mm{
R/\pi\ds R/\pi^3\ar[r]^-{\smatze{\pi}{1}}\ar[d]_{\smatze{\pi}{1}} & R/\pi^2\ar[d]^1 \\
R/\pi^2\ar[r]^1                                                   & R/\pi^2         \\
}}\right)
& &
\left(\raisebox{6mm}{\xymatrix@R-4mm{X_{21}\ar[d]_{\veps_{21}}\\\Og SX_{21}}}\right) 
& = & 
\left(\raisebox{6mm}{\xymatrix@R-4mm{
R/\pi\ar[r]\ar[d] & 0\ar[d] \\
0\ar[r]           & 0       \\
}}\right) \vsp{1} \\
\left(\raisebox{6mm}{\xymatrix@R-4mm{X_9\ar[d]_{\veps_9}\\\Og SX_9}}\right) 
& = & 
\left(\raisebox{6mm}{\xymatrix@R-4mm{
R/\pi^3\ar[r]^\pi\ar[d]_1 & R/\pi^3\ar[d]^1 \\
R/\pi^2\ar[r]^\pi         & R/\pi^3         \\
}}\right)
& &
\left(\raisebox{6mm}{\xymatrix@R-4mm{X_{22}\ar[d]_{\veps_{22}}\\\Og SX_{22}}}\right) 
& = & 
\left(\raisebox{6mm}{\xymatrix@R-4mm{
0\ar[r]\ar[d] & R/\pi\ar[d]^1 \\
0\ar[r]       & R/\pi         \\
}}\right) \\
\left(\raisebox{6mm}{\xymatrix@R-4mm{X_{10}\ar[d]_{\veps_{10}}\\\Og SX_{10}}}\right) 
& = & 
\left(\raisebox{6mm}{\xymatrix@R-4mm{
R/\pi^2\ar[r]\ar[d] & 0\ar[d] \\
0\ar[r]             & 0       \\
}}\right)
& &
\left(\raisebox{6mm}{\xymatrix@R-4mm{X_{23}\ar[d]_{\veps_{23}}\\\Og SX_{23}}}\right) 
& = & 
\left(\raisebox{6mm}{\xymatrix@R-4mm{
R/\pi^3\ar[r]^\pi\ar[d]_1 & R/\pi^2\ar[d]^1 \\
R/\pi\ar[r]^\pi           & R/\pi^2         \\
}}\right) \vsp{1} \\
\left(\raisebox{6mm}{\xymatrix@R-4mm{X_{11}\ar[d]_{\veps_{11}}\\\Og SX_{11}}}\right) 
& = & 
\left(\raisebox{6mm}{\xymatrix@R-4mm{
0\ar[r]\ar[d] & R/\pi^2\ar[d]^1 \\
0\ar[r]       & R/\pi^2         \\
}}\right)
& &
\left(\raisebox{6mm}{\xymatrix@R-4mm{X_{24}\ar[d]_{\veps_{24}}\\\Og SX_{24}}}\right) 
& = & 
\left(\raisebox{6mm}{\xymatrix@R-4mm{
R/\pi^2\ar[r]^{\pi^2}\ar[d]_1 & R/\pi^3\ar[d]^1 \\
R/\pi\ar[r]^{\pi^2}           & R/\pi^3         \\
}}\right) \\
\left(\raisebox{6mm}{\xymatrix@R-4mm{X_{12}\ar[d]_{\veps_{12}}\\\Og SX_{12}}}\right) 
& = & 
\left(\raisebox{6mm}{\xymatrix@R-4mm{
R/\pi\ds R/\pi^3\ar[r]^-{\smatze{\;\pi^2}{\pi}}\ar[d]_{\smatze{\pi}{1}} & R/\pi^3\ar[d]^1 \\
R/\pi^2\ar[r]^\pi                                                       & R/\pi^3         \\
}}\right)\vsp{1}
& &
\left(\raisebox{6mm}{\xymatrix@R-4mm{X_{25}\ar[d]_{\veps_{25}}\\\Og SX_{25}}}\right) 
& = & 
\left(\raisebox{6mm}{\xymatrix@R-4mm{
R/\pi^2\ar[r]^\pi\ar[d]_1 & R/\pi^2\ar[d]^1 \\
R/\pi\ar[r]^\pi           & R/\pi^2         \\
}}\right)\;\;\; . \\
\left(\raisebox{6mm}{\xymatrix@R-4mm{X_{13}\ar[d]_{\veps_{13}}\\\Og SX_{13}}}\right) 
& = &
\left(\raisebox{6mm}{\xymatrix@R-4mm{
R/\pi^3\ar[r]^-{\smatez{1}{\pi}}\ar[d]_1 & R/\pi\ds R/\pi^3\ar[d]^{\smatzz{1}{0}{0}{1}} \\
R/\pi^2\ar[r]^-{\smatez{1}{\pi}}         & R/\pi\ds R/\pi^3                             \\
}}\right)
& & & & \\
\ea
\]
\end{footnotesize}%

\begin{Remark}
\label{RemE3_5}\rm
Keep the assumptions of Proposition~\ref{PropE2}.

We have $(\Og\0 S)^2\, Y \iso (\Og\0 S) Y$ for $Y\in\Ob\ulmodr A$. 

The unit of the adjunction $S\adj\Og$ at an $A$-module $X\lraa{f} Y$ is represented by a factorisation 
\[
\xymatrix{
X\ar[r]^f\are[d]_{\b f\;}                & Y\ar@<-0.5mm>@{=}[d] \\
\III_f \arm[r]^*+<0mm,1mm>{\scm\dot f}   & Y                    \\
}
\]
over an image $\III_f$ of the module-defining morphism $f$. 

I do not know why.
\end{Remark}

\subsection{Another example of a left adjoint}
\label{SecExLeftAdjB}

Recall that $R$ is a principal ideal domain, with a maximal ideal generated by an element $\pi\in R$.

\subsubsection{A list of indecomposables}
\label{SecExListB}

Let 
\[
B \; :=\; A/(\pi^2 a) \= (R/\pi^3)(e\lraa{a} f)/(\pi^2 a)\; .
\]

Indecomposable nonprojective $B$-modules become indecomposable nonprojective $A$-modules via restriction along the residue class map $A\lra B$. 

We list the $24$ representatives of isoclasses of indecomposable nonprojective $B$-modules in the numbering used in \S\ref{SecExListA} as follows.
\[
\ba{l}
X_1\,,\;X_2\,,\;X_3\,,\;X_5\,,\;X_6\,,\;X_7\,,\;X_8\,,\;X_9\,,\;X_{10}\,,\;X_{11}\,,\;X_{12}\,,\;X_{13}\,,\;X_{14}\,,\;X_{15}\,,\;X_{16}\,, \\
X_{17}\,,\;X_{18}\,,\;X_{19}\,,\; X_{20}\,,\;X_{21}\,,\;X_{22}\,,\;X_{23}\,,\;X_{24}\,,\;X_{25}                                             \\
\ea
\]

\subsubsection{Construction of a left adjoint}
\label{SecConstrLeftAdjB}

Our aim in this section is to computationally verify the

\begin{Proposition}
\label{PropE4}
Suppose given a prime $p\in [2,997]$. Suppose that $R = \Fu{p}[X]$ and $\pi = X$.

Then the Heller operator $\,\Og : \ulmodr B\lra\ulmodr B$ has a left adjoint.  
\end{Proposition}

We proceed analogously to \S\ref{SecExLeftAdjA}. 

We obtain
% Data: cf. polaris.math.rwth-aachen.de:magma/sebastian/check_adjoint_5 resp. /check_adjoint_5_loop for the epsilon values further below
\[
\hsp{-10}
\ba{lclclclclclclcl}
SX_1    & = & X_2           & \hsp{-4} & \Og SX_1    & = & X_{21}           & \hsp{0} & SX_{14} & = & X_2\ds X_9        & \hsp{-4} & \Og SX_{14} & = & X_{21}\ds X_{11}        \\
SX_2    & = & X_1           &          & \Og SX_2    & = & X_3              &         & SX_{15} & = & X_2               &          & \Og SX_{15} & = & X_{21}                  \\
SX_3    & = & X_1           &          & \Og SX_3    & = & X_3              &         & SX_{16} & = & X_1\ds X_{17}     &          & \Og SX_{16} & = & X_3\ds X_{25}           \\
SX_5    & = & X_{17}        &          & \Og SX_5    & = & X_{25}           &         & SX_{17} & = & X_1\ds X_2\ds X_9 &          & \Og SX_{17} & = & X_{21}\ds X_3\ds X_{11} \\
SX_6    & = & X_2\ds X_{13} &          & \Og SX_6    & = & X_{21}\ds X_{22} &         & SX_{18} & = & X_1\ds X_2\ds X_9 &          & \Og SX_{18} & = & X_{21}\ds X_3\ds X_{11} \\
SX_7    & = & X_1\ds X_9    &          & \Og SX_7    & = & X_3\ds X_{11}    &         & SX_{19} & = & X_2\ds X_9        &          & \Og SX_{19} & = & X_{21}\ds X_{11}        \\
SX_8    & = & X_1\ds X_2    &          & \Og SX_8    & = & X_{21}\ds X_3    &         & SX_{20} & = & X_1               &          & \Og SX_{20} & = & X_3                     \\
SX_9    & = & X_{17}        &          & \Og SX_9    & = & X_{25}           &         & SX_{21} & = & X_2               &          & \Og SX_{21} & = & X_{21}                  \\
SX_{10} & = & X_2           &          & \Og SX_{10} & = & X_{21}           &         & SX_{22} & = & X_{13}            &          & \Og SX_{22} & = & X_{22}                  \\
SX_{11} & = & X_9           &          & \Og SX_{11} & = & X_{11}           &         & SX_{23} & = & X_{17}            &          & \Og SX_{23} & = & X_{25}                  \\
SX_{12} & = & X_2\ds X_{17} &          & \Og SX_{12} & = & X_{21}\ds X_{25} &         & SX_{24} & = & X_2\ds X_9        &          & \Og SX_{24} & = & X_{21}\ds X_{11}        \\
SX_{13} & = & X_1\ds X_9    &          & \Og SX_{13} & = & X_3\ds X_{11}    &         & SX_{25} & = & X_{17}            &          & \Og SX_{25} & = & X_{25}                  \\
\ea
\]
and
\begin{footnotesize}
\[
\ba{lllllll}
\left(\raisebox{6mm}{\xymatrix@R-4mm{X_1\ar[d]_{\veps_1}\\\Og SX_1}}\right) 
& = & 
\left(\raisebox{6mm}{\xymatrix@R-4mm{
R/\pi\ar[r]^1\ar[d]_1 & R/\pi\ar[d] \\
R/\pi\ar[r]           & 0           \\
}}\right)
& \hsp{5} &
\left(\raisebox{6mm}{\xymatrix@R-4mm{X_{14}\ar[d]_{\veps_{14}}\\\Og SX_{14}}}\right) 
& = & 
\left(\raisebox{6mm}{\xymatrix@R-4mm{
R/\pi^3\ar[r]^{\pi^2}\ar[d]_1 & R/\pi^3\ar[d]^1 \\
R/\pi\ar[r]^0                 & R/\pi^2         \\
}}\right) \\
\left(\raisebox{6mm}{\xymatrix@R-4mm{X_2\ar[d]_{\veps_2}\\\Og SX_2}}\right) 
& = & 
\left(\raisebox{6mm}{\xymatrix@R-4mm{
R/\pi^2\ar[r]^1\ar[d]_1 & R/\pi^2\ar[d]^1 \\
R/\pi^2\ar[r]^1         & R/\pi           \\
}}\right)
& &
\left(\raisebox{6mm}{\xymatrix@R-4mm{X_{15}\ar[d]_{\veps_{15}}\\\Og SX_{15}}}\right) 
& = & 
\left(\raisebox{6mm}{\xymatrix@R-4mm{
R/\pi^3\ar[r]\ar[d]_1 & 0\ar[d] \\
R/\pi\ar[r]           & 0       \\
}}\right) \\
\left(\raisebox{6mm}{\xymatrix@R-4mm{X_3\ar[d]_{\veps_3}\\\Og SX_3}}\right) 
& = & 
\left(\raisebox{6mm}{\xymatrix@R-4mm{
R/\pi^2\ar[r]^1\ar[d]_1 & R/\pi\ar[d]^1 \\
R/\pi^2\ar[r]^1         & R/\pi         \\
}}\right)
& &
\hsp{-12}\left(\raisebox{6mm}{\xymatrix@R-4mm{X_{16}\ar[d]_{\veps_{16}}\\\Og SX_{16}}}\right) 
& \hsp{-12}= & 
\hsp{-12}\left(\raisebox{6mm}{\xymatrix@R-4mm{
R/\pi^2\ds R/\pi^3\ar[r]^{\smatzz{1}{\pi}{1}{0}}\ar[d]_{\rsmatzz{1}{1}{1}{0}} & R/\pi\ds R/\pi^3\ar[d]^{\smatzz{1}{0}{0}{1}} \\
R/\pi^2\ds R/\pi^2\ar[r]^{\smatzz{1}{0}{0}{\pi}}                              & R/\pi\ds R/\pi^2                             \\
}}\right) \vsp{1}\\
&   & 
& &
\hsp{-12}\left(\raisebox{6mm}{\xymatrix@R-4mm{X_{17}\ar[d]_{\veps_{17}}\\\Og SX_{17}}}\right) 
& \hsp{-12}= & 
\hsp{-12}\left(\raisebox{6mm}{\xymatrix@R-4mm{
R/\pi\ds R/\pi^3\ar[r]^{\smatzz{\pi}{\;\pi^2}{1}{0}}\ar[d]_{\smatzz{1}{0}{0}{1}} & R/\pi^2\ds R/\pi^3\ar[d]^{\smatzz{1}{0}{0}{1}} \\
R/\pi\ds R/\pi^2\ar[r]^{\smatzz{0}{0}{1}{0}}                                     & R/\pi\ds R/\pi^2                               \\
}}\right)\vsp{1} \\
\left(\raisebox{6mm}{\xymatrix@R-4mm{X_5\ar[d]_{\veps_5}\\\Og SX_5}}\right) 
& = & 
\left(\raisebox{6mm}{\xymatrix@R-4mm{
R/\pi^2\ar[r]^\pi\ar[d]_1 & R/\pi^3\ar[d]^1 \\
R/\pi^2\ar[r]^\pi         & R/\pi^2         \\
}}\right)
& &
\hsp{-12}\left(\raisebox{6mm}{\xymatrix@R-4mm{X_{18}\ar[d]_{\veps_{18}}\\\Og SX_{18}}}\right) 
& \hsp{-12}= & 
\hsp{-12}\left(\raisebox{6mm}{\xymatrix@R-4mm{
R/\pi\ds R/\pi^3\ar[r]^{\smatzz{0}{\;\pi^2}{1}{\pi}}\ar[d]_{\smatzz{1}{0}{0}{1}} & R/\pi\ds R/\pi^3\ar[d]^{\rsmatzz{1}{-\pi}{0}{1}} \\
R/\pi\ds R/\pi^2\ar[r]^{\smatzz{0}{0}{1}{0}}                                     & R/\pi\ds R/\pi^2                                 \\
}}\right) \\
\left(\raisebox{6mm}{\xymatrix@R-4mm{X_6\ar[d]_{\veps_6}\\\Og SX_6}}\right) 
& = & 
\left(\raisebox{6mm}{\xymatrix@R-4mm{
R/\pi\ar[r]^\pi\ar[d]_1 & R/\pi^2\ar[d]^1 \\
R/\pi\ar[r]^0           & R/\pi           \\
}}\right)
& &
\left(\raisebox{6mm}{\xymatrix@R-4mm{X_{19}\ar[d]_{\veps_{19}}\\\Og SX_{19}}}\right) 
& = & 
\left(\raisebox{6mm}{\xymatrix@R-4mm{
R/\pi\ar[r]^{\pi^2}\ar[d]_1 & R/\pi^3\ar[d]^1 \\
R/\pi\ar[r]^0               & R/\pi^2         \\
}}\right)\vsp{1} \\
\left(\raisebox{6mm}{\xymatrix@R-4mm{X_7\ar[d]_{\veps_7}\\\Og SX_7}}\right) 
& = & 
\left(\raisebox{6mm}{\xymatrix@R-4mm{
R/\pi^2\ar[r]^-{\smatez{1}{\pi}}\ar[d]_1 & R/\pi\ds R/\pi^3\ar[d]^{\rsmatzz{1}{-\pi}{0}{1}} \\
R/\pi^2\ar[r]^-{\smatez{1}{0}}           & R/\pi\ds R/\pi^2                                 \\
}}\right)
& & 
\left(\raisebox{6mm}{\xymatrix@R-4mm{X_{20}\ar[d]_{\veps_{20}}\\\Og SX_{20}}}\right) 
& = & 
\left(\raisebox{6mm}{\xymatrix@R-4mm{
R/\pi^3\ar[r]^1\ar[d]_1 & R/\pi\ar[d]^1 \\
R/\pi^2\ar[r]^1         & R/\pi         \\
}}\right) \vsp{1}\\
\ea
\]
\[
\ba{lllllll}
\left(\raisebox{6mm}{\xymatrix@R-4mm{X_8\ar[d]_{\veps_8}\\\Og SX_8}}\right) 
& = & 
\left(\raisebox{6mm}{\xymatrix@R-4mm{
R/\pi\ds R/\pi^3\ar[r]^-{\smatze{\pi}{1}}\ar[d]_{\smatzz{1}{0}{0}{1}} & R/\pi^2\ar[d]^1 \\
R/\pi\ds R/\pi^2\ar[r]^-{\smatze{0}{1}}                               & R/\pi           \\
}}\right)
& &
\left(\raisebox{6mm}{\xymatrix@R-4mm{X_{21}\ar[d]_{\veps_{21}}\\\Og SX_{21}}}\right) 
& = & 
\left(\raisebox{6mm}{\xymatrix@R-4mm{
R/\pi\ar[r]\ar[d]_1 & 0\ar[d] \\
R/\pi\ar[r]         & 0       \\
}}\right) \vsp{1} \\
\left(\raisebox{6mm}{\xymatrix@R-4mm{X_9\ar[d]_{\veps_9}\\\Og SX_9}}\right) 
& = & 
\left(\raisebox{6mm}{\xymatrix@R-4mm{
R/\pi^3\ar[r]^\pi\ar[d]_1 & R/\pi^3\ar[d]^1 \\
R/\pi^2\ar[r]^\pi         & R/\pi^2         \\
}}\right)
& &
\left(\raisebox{6mm}{\xymatrix@R-4mm{X_{22}\ar[d]_{\veps_{22}}\\\Og SX_{22}}}\right) 
& = & 
\left(\raisebox{6mm}{\xymatrix@R-4mm{
0\ar[r]\ar[d] & R/\pi\ar[d]^1 \\
0\ar[r]       & R/\pi         \\
}}\right) \\
\left(\raisebox{6mm}{\xymatrix@R-4mm{X_{10}\ar[d]_{\veps_{10}}\\\Og SX_{10}}}\right) 
& = & 
\left(\raisebox{6mm}{\xymatrix@R-4mm{
R/\pi^2\ar[r]\ar[d]_1 & 0\ar[d] \\
R/\pi\ar[r]           & 0       \\
}}\right)
& &
\left(\raisebox{6mm}{\xymatrix@R-4mm{X_{23}\ar[d]_{\veps_{23}}\\\Og SX_{23}}}\right) 
& = & 
\left(\raisebox{6mm}{\xymatrix@R-4mm{
R/\pi^3\ar[r]^\pi\ar[d]_1 & R/\pi^2\ar[d]^1 \\
R/\pi^2\ar[r]^\pi         & R/\pi^2         \\
}}\right) \vsp{1}\\
\left(\raisebox{6mm}{\xymatrix@R-4mm{X_{11}\ar[d]_{\veps_{11}}\\\Og SX_{11}}}\right) 
& = & 
\left(\raisebox{6mm}{\xymatrix@R-4mm{
0\ar[r]\ar[d] & R/\pi^2\ar[d]^1 \\
0\ar[r]       & R/\pi^2         \\
}}\right)
& &
\left(\raisebox{6mm}{\xymatrix@R-4mm{X_{24}\ar[d]_{\veps_{24}}\\\Og SX_{24}}}\right) 
& = & 
\left(\raisebox{6mm}{\xymatrix@R-4mm{
R/\pi^2\ar[r]^{\pi^2}\ar[d]_1 & R/\pi^3\ar[d]^1 \\
R/\pi\ar[r]^0                 & R/\pi^2         \\
}}\right) \\
\left(\raisebox{6mm}{\xymatrix@R-4mm{X_{12}\ar[d]_{\veps_{12}}\\\Og SX_{12}}}\right) 
& = & 
\left(\raisebox{6mm}{\xymatrix@R-4mm{
R/\pi\ds R/\pi^3\ar[r]^-{\smatze{\;\pi^2}{\pi}}\ar[d]_{\smatzz{1}{0}{0}{1}} & R/\pi^3\ar[d]^1 \\
R/\pi\ds R/\pi^2\ar[r]^-{\smatze{0}{\pi}}                                   & R/\pi^2         \\
}}\right)
& &
\left(\raisebox{6mm}{\xymatrix@R-4mm{X_{25}\ar[d]_{\veps_{25}}\\\Og SX_{25}}}\right) 
& = & 
\left(\raisebox{6mm}{\xymatrix@R-4mm{
R/\pi^2\ar[r]^\pi\ar[d]_1 & R/\pi^2\ar[d]^1 \\
R/\pi^2\ar[r]^\pi         & R/\pi^2         \\
}}\right) \;\;\; .\vsp{1} \\
\left(\raisebox{6mm}{\xymatrix@R-4mm{X_{13}\ar[d]_{\veps_{13}}\\\Og SX_{13}}}\right) 
& = & 
\left(\raisebox{6mm}{\xymatrix@R-4mm{
R/\pi^3\ar[r]^-{\smatez{1}{\pi}}\ar[d]_1 & R/\pi\ds R/\pi^3\ar[d]^{\rsmatzz{1}{-\pi}{0}{1}} \\
R/\pi^2\ar[r]^-{\smatez{1}{0}}           & R/\pi\ds R/\pi^2                                 \\
}}\right)
& & & & \\
\ea
\]
\end{footnotesize}%

Cf.\ \S\ref{SecExLeftAdjA}.

\begin{Remark}
\label{RemE5}\rm
Keep the assumptions of Proposition~\ref{PropE4}.

We have $(\Og\0 S)^2\, Y \iso (\Og\0 S) Y$ for $Y\in\Ob\ulmodr B$. 

Given an $R/\pi^3$\nbd-module $X$, we write $\b X := X/\pi^2 X$ and
$\Ann_\pi\b X := \{\,\b x\in\b X\,:\,\pi\b x = 0\,\}$. 

The unit of the adjunction $S\adj\Og$ at a $B$-module $X\lraa{f} Y$ is represented by the composite
\[
\xymatrix@C+3mm{
X\ar[r]^f\are[d]         & Y\ar@<-0.5mm>~+{|*\dir{|}}[d]              \\
\b X\ar[r]^{\b f}\are[d] & \b Y\ar@<-0.5mm>~+{|*\dir{|}}[d]           \\
\b X/\pi\Kern\b f\ar[r]  & \b Y/(\Ann_\pi\b X)\b f\zw{,}              \\
}
\]
where the vertical maps are the respective residue class maps, and the middle and lower horizontal maps are the induced maps. 

I do not know why.
\end{Remark}

\subsection{Further examples of left adjoints}
\label{SecExLeftAdjFurther}

Let 
\[
\ba{rclcl}
C_1 & := & A/(\pi^2 f)         & = & (R/\pi^3)(e\lraa{a} f)/(\pi^2 f)       \\ % polaris.math.rwth-aachen.de:magma/sebastian/check_adjoint_2
C_2 & := & A/(\pi f)           & = & (R/\pi^3)(e\lraa{a} f)/(\pi f)         \\ % check_adjoint_3
C_3 & := & A/(\pi^2 e,\pi^2 f) & = & (R/\pi^2)(e\lraa{a} f)                 \\ % check_adjoint_4
C_4 & := & A/(\pi a)           & = & (R/\pi^3)(e\lraa{a} f)/(\pi a)         \\ % check_adjoint_6
C_5 & := & A/(\pi a,\pi^2 f)   & = & (R/\pi^3)(e\lraa{a} f)/(\pi a,\pi^2 f) \\ % check_adjoint_8
C_6 & := & A/(\pi^2 e)         & = & (R/\pi^3)(e\lraa{a} f)/(\pi^2 e)       \\ % check_adjoint_9
C_7 & := & A/(\pi e)           & = & (R/\pi^3)(e\lraa{a} f)/(\pi e)         \\ % check_adjoint_10
C_8 & := & A/(\pi^2 e,\pi a)   & = & (R/\pi^3)(e\lraa{a} f)/(\pi^2 e,\pi a) \\ % check_adjoint_11
\ea
\]

\begin{Proposition}
\label{PropE6}
Suppose given a prime $p\in [2,997]$. Suppose that $R = \Fu{p}[X]$ and $\pi = X$.

Then the Heller operator $\,\Og : \ulmodr C_j\lra\ulmodr C_j$ has a left adjoint for $j\in [1,8]$.
\end{Proposition}

\begin{Remark}
\label{RemE7}\rm
Keep the assumptions of Proposition~\ref{PropE6}.

We have $(\Og\0 S)^2\, Y \iso (\Og\0 S) Y$ for $Y\in\Ob\ulmodr C_j$ for $j\in [1,8]\ohne\{ 5\}$. 

For $j = 5$, we have 
\[
\barcl
(\Og\0 S) X_{10} & = & X_{10} \ds X_{21} \\
(\Og\0 S) X_{21} & = & X_{21}            \\
\ea
\]
in the notation of \S\ref{SecExListA}, i.e.\ 
\[
\ba{lcl}
(\Og\0 S) (R/\pi^2\lra 0) & = & (R/\pi^2\lra 0) \ds (R/\pi\lra 0) \\
(\Og\0 S) (R/\pi\lra 0)   & = & (R/\pi\lra 0)\; .                 \\
\ea
\]
\end{Remark}

\subsection{Counterexample: no right adjoint}
\label{SecExRightAdj}

Recall from \S\ref{SecExLeftAdjFurther} that $C_3 = (R/\pi^2)(e\lraa{a} f)$. As representatives of isoclasses of nonprojective $C_3$-modules we obtain,
in the notation of~\S\ref{SecExListA},
\[
\ba{lclclclclcl}
Y_1 & := & X_1    & = & (R/\pi\lraa{1} R/\pi)     & \hsp{5} & Y_5 & := & X_{21} & = & (R/\pi\lra 0)                   \\
Y_2 & := & X_3    & = & (R/\pi^2\lraa{1} R/\pi)   &         & Y_6 & := & X_{22} & = & (0\lra R/\pi)                   \\
Y_3 & := & X_6    & = & (R/\pi\lraa{\pi} R/\pi^2) &         & Y_7 & := & X_{25} & = & (R/\pi^2\lraa{\pi} R/\pi^2)\; . \\
Y_4 & := & X_{10} & = & (R/\pi\lra 0)             &         &     &    &        &   &                                 \\
\ea
\]

\begin{Remark}
\label{RemE8}
Suppose that $R = \Fu{3}[X]$ and $\pi = X$.

The functor $\,\Og : \ulmodr C_3 \lra \ulmodr C_3\,$ does not have a right adjoint.
\end{Remark}

{\it Proof.} {\sc Magma} yields
\[
H \;:=\; \big(\dim_{\Fu{3}}(\liu{\ulmodr C_3}{(Y_i\,,Y_j)})\,\big)_{i,j} \=
\enger{\left(\ba{ccccccc}
\scm 1 &\scm 0 &\scm 1 &\scm 0 &\scm 1 &\scm 0 &\scm 0 \\
\scm 1 &\scm 1 &\scm 1 &\scm 1 &\scm 1 &\scm 0 &\scm 1 \\
\scm 0 &\scm 1 &\scm 1 &\scm 0 &\scm 1 &\scm 1 &\scm 0 \\
\scm 0 &\scm 1 &\scm 0 &\scm 2 &\scm 1 &\scm 0 &\scm 1 \\
\scm 0 &\scm 1 &\scm 0 &\scm 1 &\scm 1 &\scm 0 &\scm 1 \\
\scm 1 &\scm 1 &\scm 0 &\scm 0 &\scm 0 &\scm 1 &\scm 0 \\
\scm 0 &\scm 1 &\scm 1 &\scm 1 &\scm 1 &\scm 1 &\scm 1 \\
\ea\right)}
\;\in\;(\Z_{\ge 0})^{7\ti 7} 
\]
and
\[
H' \;:=\; \big(\dim_{\Fu{3}}(\liu{\ulmodr C_3}{(\Og Y_i\,,Y_j)})\,\big)_{i,j} \=
\enger{\left(\ba{ccccccc}
\scm 1 &\scm 0 &\scm 1 &\scm 0 &\scm 1 &\scm 0 &\scm 0 \\
\scm 1 &\scm 1 &\scm 0 &\scm 0 &\scm 0 &\scm 1 &\scm 0 \\
\scm 0 &\scm 1 &\scm 1 &\scm 0 &\scm 1 &\scm 1 &\scm 0 \\
\scm 0 &\scm 0 &\scm 0 &\scm 0 &\scm 0 &\scm 0 &\scm 0 \\
\scm 0 &\scm 1 &\scm 1 &\scm 0 &\scm 1 &\scm 1 &\scm 0 \\
\scm 1 &\scm 1 &\scm 0 &\scm 0 &\scm 0 &\scm 1 &\scm 0 \\
\scm 0 &\scm 0 &\scm 0 &\scm 0 &\scm 0 &\scm 0 &\scm 0 \\
\ea\right)}
\;\in\;(\Z_{\ge 0})^{7\ti 7} \; .
\]

{\it Assume} that $\Omega$ has right adjoint $T : \ulmodr C_3\lra \ulmodr C_3\,$. 

Write $TY_j \iso \Ds_{k\in [1,7]} Y_k^{\ds u_{k,j}}$ for $j\in [1,7]$, where $U := (u_{k,j})_{k,j}\in (\Z_{\ge 0})^{7\ti 7}$. We obtain
\[
\barcl
H' 
& = & \big(\dim_{\Fu{3}}(\liu{\ulmodr C_3}{(\Og Y_i\,,Y_j)})\,\big)_{i,j}                             \\
& = & \big(\dim_{\Fu{3}}(\liu{\ulmodr C_3}{(Y_i\,, TY_j)})\,\big)_{i,j}                               \\
& = & \big(\dim_{\Fu{3}}(\liu{\ulmodr C_3}{(Y_i\,, \Ds_{k\in [1,7]} Y_k^{\ds u_{k,j}})})\,\big)_{i,j} \\
& = & \big(\sum_{k\in [1,7]} \dim_{\Fu{3}}(\liu{\ulmodr C_3}{(Y_i\,, Y_k)})\cdot u_{k,j}\big)_{i,j}   \\
& = & H \cdot U\; .                                                                                   \\
\ea
\]
So every column of $H'$ is a linear combination of columns in $H$ with coefficients in $\Z_{\ge 0}\,$. However, the third column of $H'$ would afford a coefficient $\in\Z_{> 0}$ at
the first, third or fifth column of $H$ because its first entry equals $1$. But then its second entry would also be in $\Z_{> 0}\,$, because these columns of $H$ all have second entry equal to $1$.
But this second entry equals $0$. We have arrived at a {\it contradiction.} \qed

  % Example

\parskip0.0ex
\begin{footnotesize}

\parskip1.2ex

\vspace*{3cm}

\hfill
\begin{minipage}{8cm}
\begin{flushright}
Matthias K\"unzer\\
Lehrstuhl D f\"ur Mathematik\\
RWTH Aachen\\
Templergraben 64\\
D-52062 Aachen \\
kuenzer@math.rwth-aachen.de \\
www.math.rwth-aachen.de/$\sim$kuenzer\\
\end{flushright}
\end{minipage}
\end{footnotesize}
\end{document}